\documentclass{ifacconf}

\usepackage{graphicx}      % include this line if your document contains figures
\usepackage{natbib}        % required for bibliography
\usepackage{amsmath,amssymb,amsfonts}%,bookman}
\usepackage{epsfig}
\usepackage{color}
\usepackage{caption}
\def\bm#1{\mathbf{#1}}

\newcommand{\R}{\mathbb R}

\newcommand{\Rn}{\R^n}

\newcommand{\mat}[1]{\left[\begin{array}{#1}}
\newcommand{\rix}{\end{array}\right]}
\renewcommand{\t}{^{\mbox{\tiny\sf T}}}

\newcommand{\diag}{\mathrm{diag}}
\newcommand{\sgn}{\mathrm{sgn}}

\begin{document}
\begin{frontmatter}

\title{Convergence results for continuous-time dynamics arising in ant colony optimization\thanksref{footnoteinfo}} 
% Title, preferably not more than 10 words.

\thanks[footnoteinfo]{This research was supported by CNPq (Brazil), Dept.~of Science and Technology (India), Inria (France)}

\author[First]{P.-A. Bliman}
\author[Second]{A. Bhaya}
\author[Second]{E. Kaszkurewicz} 
\author[Third]{Jayadeva} 

\address[First]{Institut National de Recherche en Informatique et en Automatique, Inria, France (e-mail: pierre-alexandre.bliman@inria.fr)}
\address[Second]{Electrical Engineering Department, Federal University of Rio de Janeiro, Brazil (e-mail: \{amit,eugenius\}@nacad.ufrj.br).}
\address[Third]{Electrical Engineering Department, Indian Institute of Technology, Delhi, India (e-mail: jayadeva@ee.iitd.ac.in)}

%\address[Fourth]{Electrical Engineering Department, Federal University of Rio de Janeiro (e-mail: amit@nacad.ufrj.br).}

\begin{abstract}                % Abstract of not more than 250 words.
This paper studies the asymptotic behavior of several continuous-time dynamical systems which are analogs of ant colony optimization algorithms that solve shortest path problems. Local asymptotic stability of the equilibrium corresponding to the shortest path is shown under mild assumptions. A complete study is given for a recently proposed model called EigenAnt: global asymptotic stability is shown, and the speed of convergence is calculated explicitly and shown to be proportional to the difference between the reciprocals of the second shortest and the shortest paths.
%The dynamics that underlies several existing ant colony optimization algorithms is elucidated in a manner that leads to new dynamics that exhibit faster convergence. For the continuous-time dynamics proposed in this paper, proofs of several important properties are given. Results of a discrete-time synchronous and asynchronous implementations for the model two-node shortest path problem illustrate the applicability of the theoretical results.
\end{abstract}


\end{frontmatter}
%===============================================================================

\section{Introduction}
Ant Colony Optimization has generated a lot of interest due to emergent optimizing behavior resulting from agents interacting through their environment, with minimal use of global information. In the most basic application of Ant Colony Optimization (ACO), a set of artificial ants find the shortest path between a source and a destination. Ants deposit pheromone on paths they take, preferring paths that have more pheromone on them. Since shorter paths are traversed faster, more pheromone accumulates on them in a given time, attracting more ants and leading to reinforcement of the pheromone trail on shorter paths. This is a positive feedback process, that can also cause trails to persist on longer paths, even though a shorter path has been later discovered and trailed by the ant colony. In \cite{ours1} and \cite{ours2}, it was shown that pheromone bias can be overcome only up to a theoretical limit; beyond that, the problem of persistence persists. ACO algorithms have employed a number of strategies to overcome this lack of plasticity. For example, in \cite{hoos} an upper bound on the amount of pheromone on a path was imposed. Finding the  optimal thresholds or parameter values is hard \citep{YuaMonBirStu2012}, and it is possible for several sub-optimal paths to end up with the maximum allowed pheromone concentration, preventing convergence to the optimal path. A more common remedial measure employed by most ACO algorithms is uniform evaporation on all paths \citep{dorigocol}. In the presence of evaporation, maintaining a trail requires continued deposition. Since evaporation necessitates the reinforcement of positive pheromone, it raises the initial bias level, on sub-optimal paths, that can be reverted by quicker returns on a shorter path. At present, it is almost ubiquitously used in most applications \citep{bulln,dorigogam,datamining,gridshop}. In the literature, opinion exists that evaporation is too slow to play an important role in foraging among real ants \citep{dorigofirstbook}; however, it is known to significantly improve the performance of artificial ant algorithms \citep{dorigobook,deneu}.

There is a large literature on ACO and its applications (see \cite{dorigobook} and references therein), but relatively less literature on the mathematical properties of ACO algorithms \citep{dorigo:theory,dorigo:survey,Blum2005}.\/ \cite{Gutjahr06} proposed a limiting process to derive a continuous-time (deterministic) differential equation from the ensemble behavior of the stochastic ant system.

One of the first ACO algorithms to provide an analysis of equilibrium states was EigenAnt, proposed in \cite{JSBKC13}. EigenAnt showed the local stability of the equilibrium corresponding to the shortest path, and presented simulation results indicating robustness of this stability to parameter choices. The approach of \cite{Gutjahr06} was used in \cite{IP12} to derive a continuous model, while EigenAnt, which is a discrete model, was proposed independently in \cite{JSBKC13}, based on similar considerations.

The robust stability properties of the EigenAnt dynamics presented in \cite{JSBKC13} motivate the question of existence of other dynamics that could display similar, or more interesting, behavior. In this context, the present paper studies continuous-time generalizations of the EigenAnt dynamics proposed in \cite{JSBKC13}, as well as the ``binary chain'' dynamics proposed in \cite{IP12}, establishing several theoretical results that were not established formally in the cited papers, notably global convergence as well as speed of convergence results. It is established herein that continuous-time EigenAnt dynamics converge globally to an equilibrium corresponding to the shortest path, i.e. from initial states corresponding to arbitrary initial pheromone concentrations on a set of paths to an equilibrium state in which all the pheromone is concentrated on the shortest path. To the best of our knowledge, this is the first such proof of global stability of the continuous-time ensemble behavior (= ODE) of any ACO algorithm. It is also important to emphasize that this global convergence is shown to be robust, in the sense that it does not depend on choices of two parameters (deposition and evaporation rates) of the algorithm.

This paper is organized as follows. The models proposed are presented in section \ref{sec_models}. A general local stability result is derived in section \ref{sec_locstab}. The particular case of the EigenAnt model is studied in section \ref{sec_eigconv}, for which the asymptotic behavior is completely described, leading, in particular, to a global stability proof and an estimate of the speed of convergence. Simulations in section \ref{sec_simul} show that some of the new variants exhibit faster convergence, demonstrating promise for use in ACO algorithms. 
%%%%%%%%%%%%%%%%%%%%%%%%%%%%%%%%%%%%%%%%%%%%%%%%%%%%%%%%%%%%%%%%%%%%%%%%%%%%%%%%%%%%%%%%%%%%%%%%%%%%
\section{EigenAnt dynamics and its variants}\label{sec_models}
The analogue of the discrete-time EigenAnt dynamics proposed in \cite{JSBKC13} is defined in continuous-time, on the positive orthant $\R_{+}^n$, as follows:
\begin{equation}\label{eq_eigant}
	\dot{\bm{x}} = \gamma\left(-\alpha \bm{I} + \frac{\beta}{\sum_i {x}_i }\bm{D}\right)\bm{x},
\end{equation}
where $\alpha, \beta, \gamma$ are scalar positive constants, $\bm{x} = (x_1 ,\ldots, x_n)$ in $\Rn$, $\bm{D} = \diag (d_1 , \ldots , d_n)$ a diagonal matrix, with diagonal entries satisfying:
\begin{equation}\label{eq_d_ordering}
	d_1 \geq d_2 \geq \cdots \geq d_n > 0
\end{equation}
In the ant colony optimization context, the $d_i$s are reciprocals of path lengths $L_i, i = 1, \ldots , n$, where $n$ paths connect a source node to a destination node and $x_i$ represents pheromone concentration on path $i$. 	In other words, (\ref{eq_eigant}) is interpreted by saying that ants deposit pheromone on path $i$, at rate $\beta d_i (x_i/\sum_i x_i)$, and it evaporates at rate $\alpha$. The number $d_1$ is the reciprocal of the shortest path length and, when the state trajectory $\bm{x}(t)$ converges to a multiple of the vector $\bm{e}_1$, this indicates that the pheromone is totally concentrated on the shortest path: in other words, the ants have `found' the shortest path. For more details on this model in the ACO context, see \cite{JSBKC13}.

%
%\color{red}
The following class of models, that generalizes the EigenAnt model, is considered in the present paper:
\begin{equation}\label{eq_phiant}
	\dot{\bm{x}} = g\left(-\alpha \bm{I} + \beta\phi(\bm{x})\bm{D}\right)\bm{x},
\end{equation}
where $\phi: \Rn \rightarrow \R$ is a real-valued differentiable function subject to the following assumptions:
\begin{itemize}
	\item[A1] As $\bm{x} \rightarrow \bm{0}$, $\phi(\bm{x}) \rightarrow \infty$.
	\item[A2] $\phi (\bm{x})$ is nonincreasing with respect to each component of $\bm{x}$.
	\item[A3] $\phi(\bm{x}) \rightarrow 0$ as $\|\bm{x}\| \rightarrow \infty$.
\end{itemize}
and $g$ is a scalar increasing differentiable function such that $g(0)=0$, and for which we define, for any diagonal matrix $\diag\{ a_1,\ldots, a_n\}$,
\begin{equation}
\label{conv}
g(\diag\{ a_1,\ldots, a_n\}) := \diag\{ g(a_1),\ldots, g(a_n)\}.
\end{equation}
Note that the choice $\phi(\bm{x}) := 1/(\sum_i x_i)$ satisfies the assumptions A1,A2, A3 and, with $g(x)= \gamma x$, converts (\ref{eq_phiant}) into the EigenAnt dynamics. As regards assumption A1, note that the choice of $\phi$ as the sum, that defines the EigenAnt dynamics, is such that $\phi(\bm{0})$ is not defined, thus this system is only studied in $\Rn\backslash\{\bm{0}\}$.
The following neural network-like variant, which use the hyperbolic tangent function and the scalar gain $\gamma>0$, and corresponds to $g=\gamma\tanh$ is also studied below:
\begin{equation}\label{eq_tanh_before}
	\dot{\bm{x}} = \gamma \left(\tanh\left(-\alpha \bm{I} + \beta\phi(\bm{x})\bm{D}\right)\right)\bm{x}
\end{equation}
%

%%%%%%%%%%%%%%%%%%%%%%%%%%%%%%%%%%%%%%%%%%%%%%%%%%%%%%%%%%%%%%%%%%%%%%%%%%%%%%%%%%%%%%%%%%%%%%%%%%%%
\section{Local asymptotic stability of the equilibria}\label{sec_locstab}
%%%%%%%%%%%%%%%%%%%%%%%%%%%%%%%%%%%%%%%%%%%%%%%%%%%%%%%%%%%%%%%%%%%%%%%%%%%%%%%%%%%%%%%%%%%%%%%%%%%%
We study the local stability behavior of the generalized model (\ref{eq_phiant}) in $\Rn\backslash\{\bm{0}\}$, in the special case when $d_i \neq d_j$ for $i \neq j$. First, the equilibria of the generalized model (\ref{eq_phiant}) are described.

It is easy to see that, under the assumptions on $\phi$ and $g$ made in Section \ref{sec_models},
 (\ref{eq_phiant}) admits exactly $n$ nonzero equilibrium points, $\bm{x}_i^*$ which, denoting the $i$th canonical vector in $\Rn$ by $\bm{e}_i$, are, for $i = 1, \ldots , n$, uniquely given by:
\begin{subequations}
\label{eq_equil_phi}
\begin{gather}
\label{eq_equil_formula}
\bm{x}_i^* = \mu_i \bm{e}_i,~\text{where $\mu_i$ is such that}\\
\label{eq_phimuei}
\phi(\mu_i \bm{e}_i) = \frac{\alpha}{\beta d_i}
\end{gather}
\end{subequations}
%
%Since $\tanh(0)=0$, (\ref{eq_tanh_before}) admits the same $n$ equilibria given in (\ref{eq_equil_formula}). In addition to these $n$ equilibria, (\ref{eq_tanh_before}) also admits the origin as an equilibrium point. Similar statements hold for (\ref{eq_sign_before}), if we assume that $\sgn(0)=0$.
For the specific choice of $\phi$ as the sum function, and $g(x)=x$, the left hand side of (\ref{eq_phimuei}) evaluates to $1/\mu_i$, yielding the explicit solution for the equilibria as $\frac{\beta d_i}{\alpha}\bm{e}_i, i = 1, \ldots , n$, as in the discrete-time case studied in~\cite{JSBKC13}.

In regard to the equilibrium point analysis just carried out, it should be emphasized that the proposed models all have the desired equilibrium point $\bm{x}_1^*$ (which corresponds to the shortest path) as one possible equilibrium, amongst others. The stability analysis, to be presented in what follows, shows that only the desired equilibrium is stable, while all others are unstable. Curiously, such analyses are virtually absent from the ACO literature, exceptions being the papers of \cite{JSBKC13,IP12}. In fact, convergence to spurious equilibria is often reported in the literature and a specific analysis of this is given in a particular case in \cite{IP12}, in which the ACO dynamics actually possess spurious stable equilibria.

The stability properties of the equilibria are derived in the following theorem.
%
%\color{red}
\begin{thm}\label{thm_equils}
Assume that $g$ is a scalar increasing differentiable function such that $g(0)=0$ and that $\phi$ satisfies the assumptions A1 through A3 given in Section 2.
Then
\begin{itemize}
\item
The equilibrium points of \eqref{eq_phiant} are exactly the $n$ distinct solutions $\bm{x}_i^* = \mu_i \bm{e}_i, i = 1, \ldots , n$ of \eqref{eq_equil_phi}.
\item
If $d_i \neq d_j$ for $i \neq j$, the equilibrium point $\bm{x}_1^*$ is locally asymptotically stable provided that $\frac{\partial\phi}{\partial x_1}(\bm{x}_1^*)<0$, while all other equilibria $\bm{x}_i^*$, $i=2,\dots, n$, are unstable.
\end{itemize}
\end{thm}

\begin{pf}
The proof is based on the analysis of the eigenvalues of the linearized system.

For the sake of simplicity, we first consider the case where $g(x)\equiv x$.
In this case, \eqref{eq_phiant} simplifies to
\begin{equation}\label{eq_H_def}
\dot{\bm{x}}
= \left(
-\alpha\bm{I} + \beta\phi(\bm{x})\bm{D}
\right) \bm{x} =: H(\bm{x}).
\end{equation}
The Jacobian of $H$ is given by:
\begin{equation}
\label{eq_Jac}
\nabla H(\bm{x}) = -\alpha\bm{I} + \beta\phi(\bm{x})\bm{D}
+\beta D\bm{x}\cdot \nabla\phi(\bm{x})\t
\end{equation}
where $\nabla\phi(\bm{x}) := \begin{pmatrix}
\frac{\partial\phi}{\partial x_1} & \dots & \frac{\partial\phi}{\partial x_n}
\end{pmatrix}\t$.
Let us study the value of $\nabla H(\bm{x}_i^*)$.
On the one hand, due to \eqref{eq_equil_phi}, the matrix $-\alpha\bm{I} + \beta\phi(\bm{x}_i^*)\bm{D}$ is a diagonal matrix whose $(j,j)$ element is $\alpha (\frac{d_j}{d_i}-1)$, $j=1,\dots, n$.
The latter are thus nonnegative for $j=1,\dots, i-1$, zero for $j=i$ and nonpositive for $j=i+1,\dots, n$.
On the other hand, the matrix $\beta D\bm{x}_i^*\cdot \nabla\phi(\bm{x}_i^*)\t$ is null, except for the $i$-th row which contains nonpositive elements (due to assumption A2).

The lower block-triangular structure of the matrix $\nabla H(\bm{x}_i^*)$ is now evident.
When $i>1$, the diagonal top-left block contains positive values, which implies instability of the corresponding equilibrium points $\bm{x}_i^*$, $i=2,\dots, n$. For $i=1$, the matrix $\nabla H(\bm{x}_i^*)$ is upper triangular with the diagonal elements $\frac{\partial\phi}{\partial x_1}(\bm{x}_1^*), \alpha (\frac{d_2}{d_1}-1),\dots, \alpha (\frac{d_n}{d_1}-1)$, and its spectrum is thus located in the open left-half complex plane whenever $\frac{\partial\phi}{\partial x_1}(\bm{x}_1^*)<0$. This shows the local asymptotic stability of $\bm{x}_1^*$ in the case where $g(x)\equiv x$.

For general functions $g$ satisfying the hypotheses of the theorem statement, the expression of the Jacobian provided in formula \eqref{eq_Jac} has to be replaced by
\[
g'(H(\bm{x}))
\left(
-\alpha\bm{I} + \beta\phi(\bm{x})\bm{D}
+\beta D\bm{x}\cdot \nabla\phi(\bm{x})\t
\right)
\]
where $H(\bm{x})$ is defined in (\ref{eq_H_def}), and, in agreement with our convention defined in \eqref{conv}, $g'(H(\bm{x})) = g'(-\alpha\bm{x} + \beta\phi(\bm{x})\bm{Dx})$ is the diagonal matrix whose elements are obtained by applying $g'$ to the diagonal elements of the diagonal matrix $-\alpha\bm{x} + \beta\phi(\bm{x})\bm{Dx}$.
The analysis is thus conducted as before, using the positivity of these coefficients, due to the fact that $g$ is an increasing function. \qed
\end{pf}

\subsection{Phase portraits of the EigenAnt and MaxAnt dynamics for $n=2$}
%%%
For the purposes of comparison with the phase plane portraits presented in \cite{IP12}, we will show the phase portraits of EigenAnt, as well as of \eqref{eq_phiant} with the choice $\phi = 1/\max\{\bm{x}\}$ (which we call the MaxAnt dynamics). The MaxAnt dynamics sets a theoretical limit on speed of convergence, and is shown here for comparison, even though it obviously cannot be used legitimately in the shortest path problem.
The notable feature, common to all the dynamics proposed in this paper, is the absence of spurious equilibria, which occur in \cite{IP12}, in cases where a certain exponent, called the pheromone amplification factor, differs from unity.
\begin{center}
	\includegraphics[width=8.4cm]{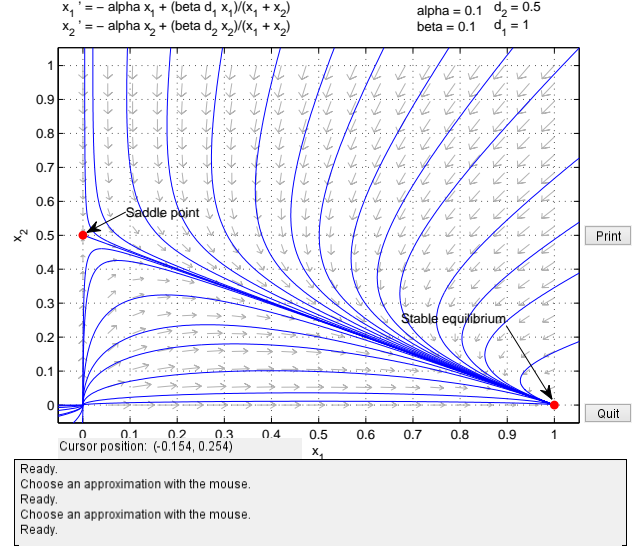}
	\captionof{figure}{Phase portrait of EigenAnt dynamics}
	\label{fig:EA_pp}
\end{center}
\begin{center}
	\includegraphics[width=8.4cm]{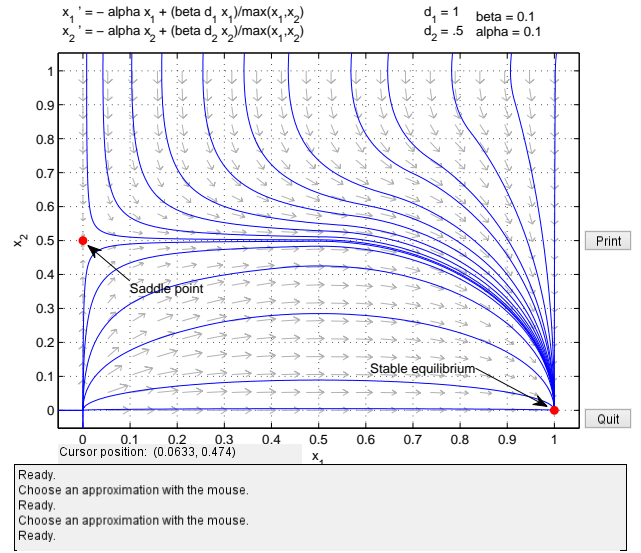}
	\captionof{figure}{Phase portrait of MaxAnt dynamics}
	\label{fig:MA_pp}
\end{center}
%%%%%%%%%%%%%%%%%%%%%%%%%%%%%%%%%%%%%%%%%%%%%%%%%%%%%%%%%%%%%%%%%%%%%%%%%%%%%%%%%%%%%%%%%%%%%%%%%%%%%
\section{Convergence properties of continuous-time EigenAnt dynamics}\label{sec_eigconv}
This section is devoted to a complete analysis of the asymptotic behavior of the EigenAnt dynamics (\ref{eq_eigant}).
It is clear that any component of the state $\bm{x}$ which departs initially from zero remains zero at any time.
Thus, with no loss of generality, we assume positive initial conditions, i.e.:
\begin{equation}
\label{eq3}
\forall i=1,\dots, n,\quad x_i(0) >0.
\end{equation}
The results given below (Theorems \ref{thm_globstab} and \ref{thm_convergence}) extend the results of Theorem \ref{thm_equils}, for EigenAnt dynamics. Several technical results (Propositions \ref{p1} to \ref{p3}) are needed in order to prove the main global stability result, Theorem \ref{thm_globstab}.
\begin{prop} [Invariance of the positive orthant]
\label{p1}
For any $i=1,\dots, n$, for any $t \geq 0$, $x_i(t)>0$.

In particular, defining
\begin{equation}
\label{eq4}
S(t) := \sum_i x_i(t)\ ,
\end{equation}
one has $S(t) >0$ for any $t \geq 0$.
\end{prop}
%Notice that, when $\alpha>1$, this property was not systematically verified by the discrete-time system.

\begin{pf}
The components of the solutions of \eqref{eq_eigant} are continuous with respect to time, and start from positive values.
As long as every component is positive, one has
\begin{equation}
\dot x_i + \alpha x_i = \frac{\beta d_i}{\sum_j x_j} x_i >0\ ,
\end{equation}
whence:
\begin{equation}
x_i(t) \geq x_i(0) e^{-\alpha t} >0
\end{equation}
and the conclusion holds. \qed \end{pf}

\begin{prop} [Upper and lower bounds of sum of states]
The following bounds hold:
\begin{equation}
\label{eq16}
\frac{\beta d_n}{\alpha} \leq \liminf_{t\to +\infty} S(t) \leq \limsup_{t\to +\infty} S(t) \leq \frac{\beta d_1}{\alpha}\ .
\end{equation}
In particular, $S$ is uniformly bounded from above, and
\begin{subequations}
\begin{gather}
\label{eq5a}
\frac{1}{S} \text{ is locally integrable on } \R^+\\
\label{eq5b}
\lim_{t\to +\infty} \int_0^t \frac{ds}{S(s)} = +\infty
\end{gather}
\end{subequations}
\end{prop}

\begin{pf}
Summing the $n$ expressions in \eqref{eq_eigant} over $i$ yields
\begin{equation}
\dot S = -\alpha S + \beta \sum_i d_i\frac{x_i}{S}\ .
\end{equation}
Thus
\begin{equation}
\beta d_n \leq \dot S + \alpha S \leq \beta d_1\ ,
\end{equation}
and therefore by integration:
\begin{align}
\frac{\beta d_n}{\alpha} + \left(
S(0) - \frac{\beta d_n}{\alpha}
\right) e^{-\alpha t}
\leq
S(t)  \nonumber \\
\leq
\frac{\beta d_1}{\alpha} + \left(
S(0) - \frac{\beta d_1}{\alpha}
\right) e^{-\alpha t} &
\end{align}
from which \eqref{eq16} is deduced.
Function $S(t)$, being positive and continuous, takes on values bounded away from zero on any compact set of $\R^+$.
Due to \eqref{eq16}, it is bounded away from zero on the whole set $\R^+$, and this in particular yields \eqref{eq5a}.
Finally, identity \eqref{eq5b} is deduced from the fact that $S$ is uniformly bounded from above. \qed \end{pf}

The following technical propositions are needed to establish the first theorem on asymptotic properties of the trajectories.
\begin{prop}
\label{p3}
For any trajectory of system \eqref{eq_eigant}, define\footnote{Note that $F$ depends on the $x_i(0)$. This dependence is not made explicit, for notational simplicity, since no confusion should arise.}
\begin{equation}
\label{eq7}
F(x) := \sum_i \frac{x_i(0)}{\beta d_i} e^{\beta d_i x}\ .
\end{equation}
Then,
\begin{itemize}
\item $F: [0,+\infty)\to [\sum_i \frac{x_i(0)}{\beta d_i},+\infty)$ is increasing and invertible;
\item for any $t>0$,
\begin{subequations}
\begin{gather}
\label{eq6a}
\int_0^t \frac{ds}{S(s)} = F^{-1} \left(
F(0) + \frac{1}{\alpha}(e^{\alpha t}-1)
\right)\\
\label{eq6b}
S(t) = e^{-\alpha t} F'\left(
\int_0^t \frac{ds}{S(s)}
\right)\\
\label{eq6c}
= e^{-\alpha t} F'\left(
F^{-1} \left(
F(0) + \frac{1}{\alpha}(e^{\alpha t}-1)
\right) \right)
\end{gather}
\end{subequations}
\item
$S(t)$ admits a (positive and finite) limit for $t\to +\infty$.
\end{itemize}
\end{prop}

\begin{pf}
Instead of summing up the $n$ identities in \eqref{eq_eigant} as was done before, we first integrate them, to obtain
\begin{equation}
\label{eq10}
x_i(t) = x_i(0) e^{\int_0^t \left(
\frac{\beta d_i}{S(s)}-\alpha
\right)\ ds}\ .
\end{equation}
Summing up now gives
\begin{equation}
S(t) = \sum_i x_i(t) = e^{-\alpha t} \sum_i x_i(0) e^{\beta d_i\int_0^t \frac{ds}{S(s)}}\ ,
\end{equation}
that is
\[
S(t) = e^{-\alpha t} F'\left(
\int_0^t \frac{ds}{S(s)}
\right)\ ,
\]
which is the first identity in \eqref{eq6b}.
One deduces
\begin{equation}
e^{\alpha t} = \frac{1}{S(t)}
F'\left(
\int_0^t \frac{ds}{S(s)}
\right)\ ,
\end{equation}
and by integration over $[0,t]$, $t>0$,
\begin{equation}
F\left(
\int_0^t \frac{ds}{S(s)}
\right) - F(0) = \frac{1}{\alpha}(e^{\alpha t} -1)\ .
\end{equation}
This gives \eqref{eq6a}, and subsequently \eqref{eq6c}.

Observing equation \eqref{eq6c} shows that $S(t)$ admits a limit for $t\to +\infty$, and this finishes the proof. \qed \end{pf}

The following theorem, which describes the overall asymptotic behavior, shows that if $d_i < d_1$, then the $i$th component of the vector $\bm{x}$ tends to zero, whereas, if $d_1 = d_2 = \cdots = d_j$, then the sum of the components $\sum_{k=1}^j x_j (t)$ tends to a fixed value as $t$ tends to infinity. Thus, in particular, it establishes global stability of the equilibrium set corresponding to the shortest path, without the assumption $d_i \neq d_j$ for $i \neq j$, thus generalizing the local stability result of Theorem \ref{thm_equils} for EigenAnt dynamics.
\begin{thm}[Global stability and asymptotic behavior]\label{thm_globstab}\mbox{}

For any solution of (\ref{eq_eigant}), for $i$ such that $d_i<d_1$:
\begin{equation}
\label{eq11}
\lim_{t\to +\infty} x_i(t) = 0\ .
\end{equation}
Moreover,
\begin{equation}
\label{eq12}
\lim_{t\to +\infty} S(t)
= \lim_{t\to +\infty} \sum_{j\ :\ d_j = d_1} x_j(t)
= \frac{\beta d_1}{\alpha}\ .
\end{equation}
\end{thm}
\begin{pf}[Theorem \ref{thm_globstab}]
From \eqref{eq6a} and using the definition of $F$ in \eqref{eq7}, one deduces
\begin{multline}
\limsup_{t\to +\infty}
e^{-\alpha t} F\left(
\int_0^t \frac{ds}{S(s)}
\right)
\\
= \limsup_{t\to +\infty}
e^{-\alpha t} \sum_i \frac{x_i(0)}{\beta d_i} e^{\beta d_i \int_0^t \frac{ds}{S(s)}}
 < +\infty\ . \label{eq8}
\end{multline}
In particular, whenever $d_i<d_1$, one has
\begin{multline}
\label{eq9}
e^{-\alpha t} e^{\beta d_i \int_0^t \frac{ds}{S(s)}}
= e^{-\beta (d_1-d_i) \int_0^t \frac{ds}{S(s)}}
\left(
e^{-\alpha t} e^{\beta d_1 \int_0^t \frac{ds}{S(s)}}
\right)\\
\leq \delta\ e^{-\beta (d_1-d_i) \int_0^t \frac{ds}{S(s)}}
\left(
e^{-\alpha t}
\sum_j \frac{x_j(0)}{\beta d_j} e^{\beta d_j \int_0^t \frac{ds}{S(s)}}
\right)
\end{multline}
where we have put
\[
\delta := \max_j \frac{\beta d_j}{x_j(0)}\ .
\]

On the other hand, from \eqref{eq5b}:
\begin{equation}
\lim_{t\to +\infty} e^{-\beta (d_1-d_i) \int_0^t \frac{ds}{S(s)}} = 0\ .
\end{equation}
Thus from \eqref{eq8} and \eqref{eq9}
\begin{align}
\lim_{t\to +\infty}
e^{-\alpha t} e^{\beta d_i \int_0^t \frac{ds}{S(s)}}
\leq \delta \lim_{t\to +\infty} e^{-\beta (d_1-d_i) \int_0^t \frac{ds}{S(s)}}\cdot \nonumber \\
\limsup_{t\to +\infty}
e^{-\alpha t} \sum_i \frac{x_i(0)}{\beta d_i} e^{\beta d_i \int_0^t \frac{ds}{S(s)}}
= 0\ . & \nonumber
\end{align}
Recognizing this expression in the right-hand side of \eqref{eq10} now yields
\[
\lim_{t\to +\infty} x_i(t) = 0
\]
which is exactly \eqref{eq11}.

The behavior at infinity of $S(t) = \sum_i x_i(t)$ is thus the same as the behavior at infinity of $\sum_{j\ :\ d_j = d_1} x_j(t)$.
Notice that, due to Proposition \ref{p3}, both functions converge.
Exploiting this fact and summing up the equations in \eqref{eq_eigant} that correspond to this smaller set of indices, it follows that, $\forall\varepsilon >0, \exists T>0, \forall t\geq T$:
\begin{multline}
\left|
\frac{d}{dt}
\left(
\sum_{j\ :\ d_j = d_1} x_j(t)
\right) \right.
+ \\
\left. \alpha \left(
\sum_{j\ :\ d_j = d_1} x_j(t)
\right)
- \beta d_1 \right|
\leq \varepsilon
\end{multline}
Removing the absolute value, integrating the differential inequalities and passing to the limit in time yields
\begin{multline}
\forall\varepsilon >0,
\frac{\beta d_1-\varepsilon}{\alpha}
\leq \liminf_{t\to +\infty} \sum_{j\ :\ d_j = d_1} x_j(t) \\
\leq \limsup_{t\to +\infty} \sum_{j\ :\ d_j = d_1} x_j(t)
\leq \frac{\beta d_1+\varepsilon}{\alpha} ,
\end{multline}
and finally \eqref{eq12} after letting $\varepsilon \to 0$.
\qed 
\end{pf}

The remainder of this section is devoted to estimating the speed of convergence to the equilibrium exhibited in Theorem \ref{thm_globstab}.
We first provide a technical result (Proposition \ref{p4}) and then state the key results in Theorem \ref{thm_convergence}, which fully describes the manner in which each component of the state vector (= pheromone concentration) tends towards its limit.

From now on, let $n'$ be the cardinality of the set $\{d_i\ :\ i=1,\dots n\}$.
We consider the positive numbers $d'_i$, $i=1,\dots n'$ such that $d_1 = d'_1 > d'_2> \dots > d'_{n'}$ and
$\{d'_i\ :\ i=1,\dots n'\} = \{d_i\ :\ i=1,\dots n\}$.
In particular (as $d'_1 = d_1$):
\begin{equation}
d'_i = \max\{ d_i\ : \ d_i < d'_{i-1} \},\quad i=1,\dots, n'
\end{equation}
\begin{prop}
\label{p4}
Denote
\begin{equation}
\label{eq18}
\sigma_i = \frac{1}{\beta d_i}\sum_{j\ :\ d_j = d'_i} x_j(0), \qquad i=1,\dots, n'
\end{equation}
With $F$ defined as in \eqref{eq7}, the following asymptotic expansions hold:
\begin{subequations}
\label{eq13}
\begin{multline}
\label{eq13a}
F^{-1}(y)
= \frac{1}{\beta d'_1} \log A ,~\mbox{where}\\
A = \left(
\frac{1}{\sigma_1} \left(
y - \sigma_2 \left(
\frac{y}{\sigma_1}
\right)^{\frac{d'_2}{d'_1}}
\left(
1+ \varepsilon(y)
\right) \right) \right)
\end{multline}
\begin{multline}
\label{eq13b}
\int_0^t \frac{ds}{S(s)}
= F^{-1}\left(
F(0) + \frac{1}{\alpha} (e^{\alpha t} -1)
\right) \\
= \frac{1}{\beta d'_1} \log B ,~\mbox{where}\\
B = \left(
\frac{1}{\sigma_1} \left(
\frac{e^{\alpha t}}{\alpha} - \sigma_2 \left(
\frac{e^{\alpha t}}{\alpha \sigma_1}
\right)^{\frac{d'_2}{d'_1}}
\left(
1+ \varepsilon(t)
\right) \right) \right)
\end{multline}
\end{subequations}
where (\ref{eq13b}) holds along any solution of (\ref{eq_eigant}) and where $\varepsilon$, in both \eqref{eq13a} and \eqref{eq13b}, denotes a function which vanishes when its argument goes to $+\infty$.
\end{prop}
\begin{pf}[Proposition \ref{p4}]
First, notice that
\begin{equation}
y \to +\infty \Leftrightarrow F^{-1}(y) \to +\infty\ .
\end{equation}
Denote $x=F^{-1}(y)$ for some $y\in\R^+$: by definition,
\begin{equation}
\label{eq15}
y = F(x) = \sum_{i=1}^n \frac{x_i(0)}{\beta d_i} e^{\beta d_i x} \ ,
\end{equation}
and thus (since $d'_1 = d_1$)
\begin{equation}
y = e^{\beta d_1 x} \sum_{i=1}^{n'} \frac{\sigma_i}{\beta d'_i} e^{\beta (d'_i-d_1)x}\ .
\end{equation}
From the fact that $d'_i<d_1$ for $i=2,\dots, n'$, one thus deduces the first order term, namely :
\begin{equation}
y = \sigma_1 e^{\beta d_1 x} (1+\varepsilon (x)) \ ,
\end{equation}
or equivalently
\begin{equation}
y (1+\varepsilon (y))= \sigma_1 e^{\beta d_1 x}\ ,
\end{equation}
or again
\begin{equation}
\label{eq14}
x = \frac{1}{\beta d_1} \log\left(
\frac{y}{\sigma_1} (1+\varepsilon (y))
\right)\ .
\end{equation}
In the formulas above, and in what follows, $\varepsilon$ denotes various functions vanishing at infinity.

Introducing \eqref{eq14} in \eqref{eq15} yields
\begin{equation}
\sigma_1 e^{\beta d_1 x} = y - \sigma_2 \left(
\frac{y}{\sigma_1}
\right)^{\frac{d'_2}{d'_1}} (1+\varepsilon(y))
\end{equation}
which gives \eqref{eq13a}.
Now putting
\[
y = F(0) +\frac{1}{\alpha} (e^{\alpha t} - 1)
\]
in \eqref{eq13a} yields \eqref{eq13b} after deletion of lower order terms.
\qed \end{pf}
\begin{thm}[Convergence rates]\label{thm_convergence}
The following expansions hold, for $t\to +\infty$, for any trajectory of the system (\ref{eq_eigant}):
\begin{subequations}
\label{eq17}
\begin{multline}
\label{eq17b}
\hspace{-.3cm}
x_i(t) = \\
x_i(0) \left(
\frac{1}{\alpha\sigma_1} - \frac{\sigma_2}{\sigma_1}\left(
\frac{1}{\alpha\sigma_1}
\right)^{\frac{d'_2}{d_1}}
e^{-\alpha \left(
1-\frac{d'_2}{d_1}
\right) t} (1+\varepsilon (t))
\right) \\ \text{ if } d'_i = d_1
\end{multline}
%\color{red}
\vspace{-1cm}
\begin{multline}
\label{eq17a}
x_i(t)
= x_i(0) \left(
\frac{1}{\alpha\sigma_1}
\right)^{\frac{d'_i}{d_1}}
e^{-\alpha \left(
1-\frac{d'_i}{d_1}
\right) t} (1+\varepsilon (t))\\ \text{ if } d'_i < d_1
\end{multline}
%\color{black}
\vspace{-1cm}
\begin{multline}
\label{eq17c}
S(t) = \sum_{i=1}^n x_i (t)\\
= \frac{\beta d_1}{\alpha}
- \beta\sigma_2 (d_1-d'_2) \left(
\frac{1}{\alpha\sigma_1}
\right)^{\frac{d'_2}{d_1}} \times \\
e^{-\alpha \left(
1-\frac{d'_2}{d_1}
\right) t} (1+\varepsilon (t))
\end{multline}
\end{subequations}
where $d_1$ and $d'_2$ are respectively the inverse of the shortest and second shortest paths (counted without multiplicity); and where the quantities $\sigma_1$ and $\sigma_2$ are defined as functions of the initial conditions in \eqref{eq18}.
\end{thm}

\begin{pf}
Formulas \eqref{eq17b} and \eqref{eq17a} are obtained by putting \eqref{eq13b} in \eqref{eq10}. The value of the sum $S(t)$ is then deduced by summing up over all indices $i$. It turns out that the terms $x_i(t)$ whose index $i$ is such that $d_i < d'_2$ can be omitted, as they lead to quantities of faster convergence. The expressions $x_i(0)$ are finally removed with the help of \eqref{eq18}, to obtain formula \eqref{eq17c}. \qed \end{pf}

Theorem \ref{thm_convergence} establishes that the components $x_i$ which do not correspond to the shortest path (i.e. $d_i < d_1$) go to zero with a decay rate proportional to $d_1-d_i$.
On the contrary, the components $x_i$ for which $d_i=d_1$ converge to the value $\frac{x_i(0)}{\alpha\sigma_1}>0$: in other words, the proportion $\frac{\lim_{t\to +\infty} x_i(t)}{\lim_{t\to +\infty} x_{i'}(t)}$ is equal to $\frac{x_i(0)}{x_{i'}(0)}$ for any paths $i,i'$ whose lengths are equal to the length of the shortest path.
Finally, the convergence occurs at a speed proportional to the difference between the shortest and the second shortest paths.
%%%%%%%%%%%%%%%%%%%%%%%%%%%%%%%%%%%%%%%%%%%%%%%%%%%%%%%%%%%%%%%%%%%%%%%%%%%%%%%%%%%%%%%%%%%%%%%%%%%
\section{Simulations of EigenAnt dynamics and its variants}\label{sec_simul}
This section presents simulations of the various models \eqref{eq_phiant} studied above. The example we use is a ten-path two-node shortest path problem, meaning that two nodes are connected by ten paths of lengths varying from $1$ to $10$. This means that diagonal matrix $\bm{D}$, which contains reciprocals of path lengths, is given by $\diag (1, 1/2,1/3, \ldots, 1/9, 1/10)$. The initial condition is chosen as $(0.1, 0.2,0.3,\ldots,0.9,1.0)$, which is referred to in ACO terminology as an initial bias (largest state or pheromone concentration on longest path).

The parameter choices are $\alpha = \beta = 1$ and $\gamma = 10$ for (\ref{eq_eigant}).
We also study model \eqref{eq_tanh_before} with $\alpha = \beta = 0.1$ and $\gamma = 10$.
We also consider the latter model with an ``infinite gain", namely
\begin{equation}\label{eq_sign_before}
	\dot{\bm{x}} = \gamma\left(\sgn\left(-\alpha \bm{I} + \beta\phi(\bm{x})\bm{D}\right)\right)\bm{x}
\end{equation}
with $\alpha = \beta = 0.1$ and $\gamma = 0.5$. Finally, for the purposes of comparison, we will also show simulations with the choice $\phi = 1/\max\{\bm{x}\}$, which sets a theoretical limit on speed of convergence, although it cannot be legitimately used in the shortest path problem.

For all simulations, integration is carried out for a horizon of $1000$ or $2000$ steps, as specified in the figure captions, using the forward Euler method with stepsize $= 0.02 $. 

%We start by comparing the EigenAnt and Maxant models in figures \ref{fig:eigenant_sync} and \ref{fig:maxant_sync}, maintaining the same values of the parameters $\alpha= \beta  = 1, \gamma=1$ and the same initial condition $(0.1, 0.2,0.3,\ldots,0.9,1.0)$. The MaxAnt system converges faster than the EigenAnt system to the stable equilibrium corresponding to the shortest path. Another notable feature is that the trajectory of the EigenAnt algorithm displays an initial decrease in the pheromone concentration ($x_1$) corresponding to the shortest path, whereas this undesirable ``nonminimum phase'' behavior is absent in the corresponding MaxAnt trajectory. 

\begin{figure}[h]
\begin{center}
\includegraphics[width=8.4cm]{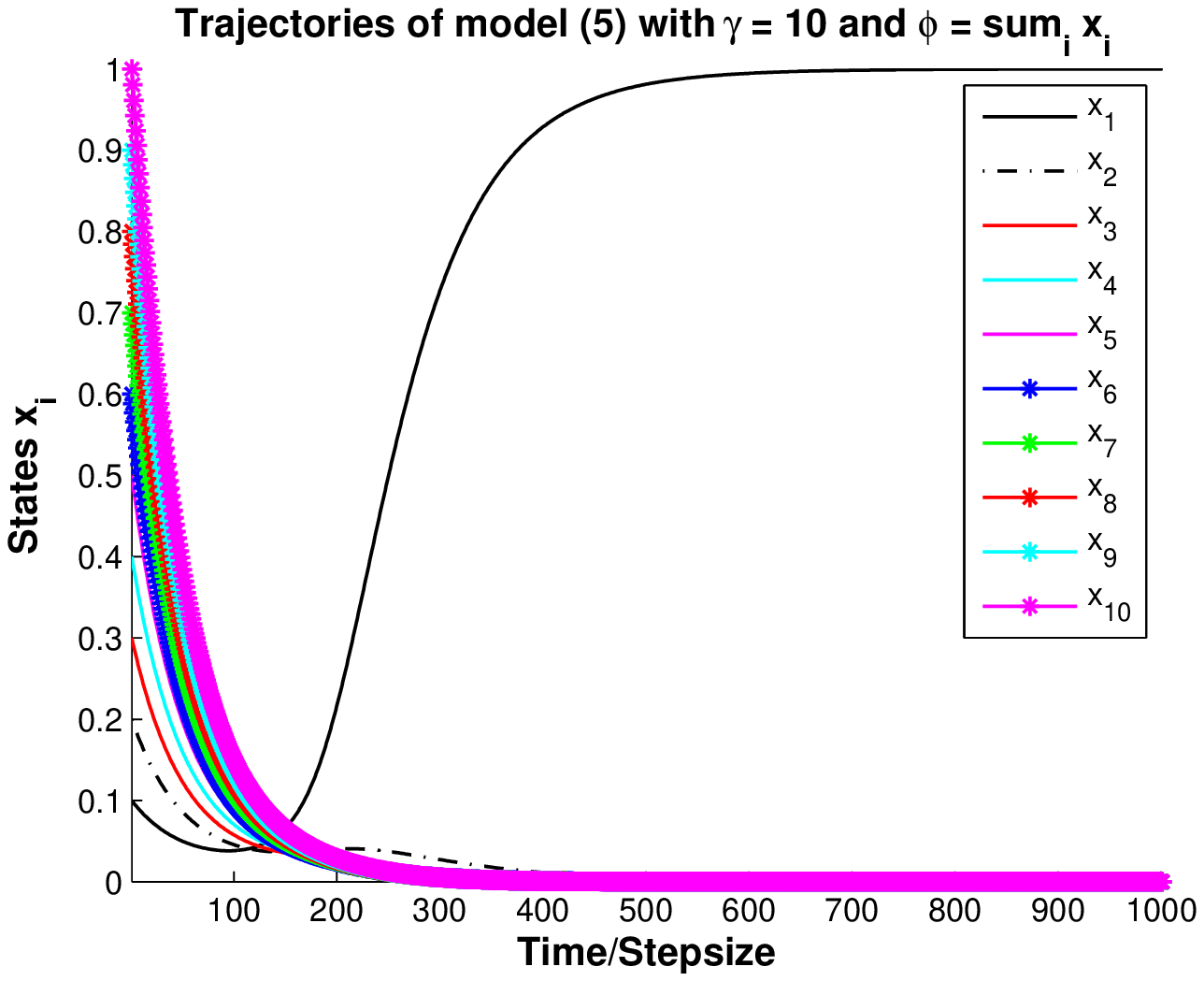}    % The printed column width is 8.4 cm.
\caption{Convergence to equilibrium corresponding to the shortest path, using hyperbolic tangent dynamics (\ref{eq_tanh_before}) with $\phi = \sum_i x_i$, for a ten-path two-node shortest path problem described at the beginning of Section \ref{sec_simul}.}
\label{fig:tanhbant_sum_sync}
\end{center}
\end{figure}

\begin{figure}
\begin{center}
\includegraphics[width=8.4cm]{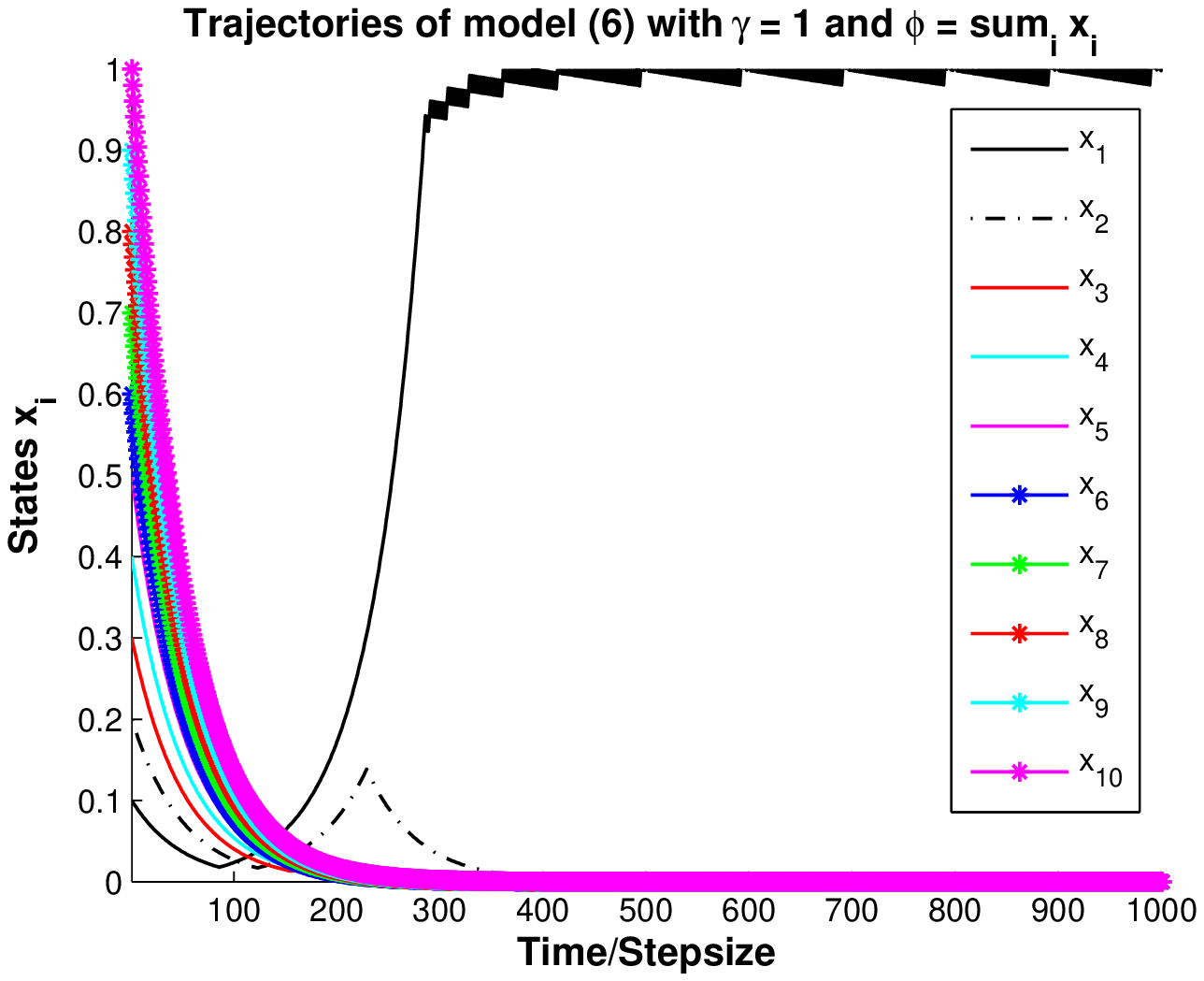}    % The printed column width is 8.4 cm.
\caption{Convergence to equilibrium corresponding to the shortest path, using signum dynamics (\ref{eq_sign_before}), with $\phi = \sum_i x_i$, for a ten-path two-node shortest path problem described at the beginning of Section \ref{sec_simul}.}
\label{fig:sigbant_sum_sync}
\end{center}
\end{figure}

Two other choices of dynamics, namely, model (\ref{eq_tanh_before}) and (\ref{eq_sign_before}), both with $\phi = \sum$, are shown in Figures \ref{fig:tanhbant_sum_sync} and figure \ref{fig:sigbant_sum_sync}, respectively. Both show fast convergence to the stable equilibrium corresponding to the shortest path. For the choice of $g = \tanh$, it should be observed that the speed of convergence depends on the parameter $\gamma$ which was chosen as $10$ in Figure \ref{fig:tanhbant_sum_sync}. On the other hand, in the case of choice of $g$ as a signum, it is necessary to set the ``external'' gain $\gamma = 1$ in order to limit the chattering, visible in Figure \ref{fig:sigbant_sum_sync}, due to numerical integration of the discontinuous right hand side of (\ref{eq_sign_before}). 

\begin{figure}[h]
\begin{center}
\includegraphics[width=8.4cm]{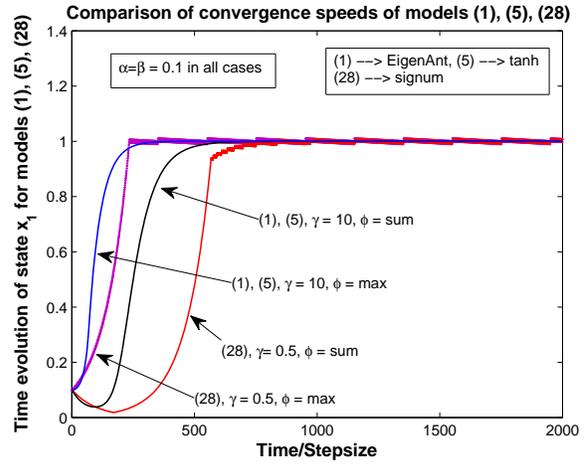}    % The printed column width is 8.4 cm.
\caption{Comparisons of speed of convergence of the state corresponding to the shortest path, for models (\ref{eq_eigant}),(\ref{eq_tanh_before}),(\ref{eq_sign_before}), for the choices $\phi = \sum_i x_i , \max\{x_i\}$.}
\label{fig:comparisons}
\end{center}
\end{figure}

Given these choices of gain $\gamma$ described above, a comparison of the trajectories of the pheromone concentration ($x_1$) corresponding to the shortest path for the six possible dynamical systems introduced above, is given in Figure \ref{fig:comparisons}. From the simulations, it appears that model (\ref{eq_phiant}) with $g(x)=x$ and with $g(x) = \tanh(x)$ have virtually indistinguishable behavior for the same choices of $\alpha,\beta,\phi$ and corresponding choices of $\gamma$. For these models, using $\phi = \max$ has a clear edge over those using $\phi = \sum$, in terms of speed of convergence and also in terms of not displaying an initial decrease in value. The model (\ref{eq_sign_before}) using the signum function can be regarded as an ``infinite'' gain version of model (\ref{eq_tanh_before}), and therefore converges the fastest of all, although the solution displays chattering. Furthermore, while it is natural to attempt to simulate \eqref{eq_sign_before} in order to speed up the convergence, the results should be interpreted cautiously in view of the fact that well-posedness of the initial value problem for equation \eqref{eq_sign_before} has not yet been established.

\begin{figure}[h]
\begin{center}
\includegraphics[width=8.4cm]{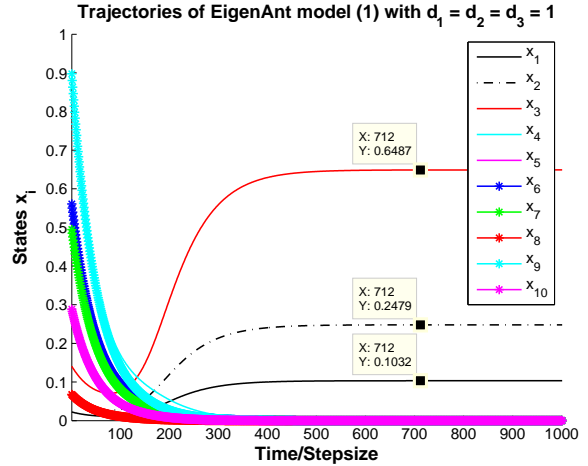}    % The printed column width is 8.4 cm.
\caption{Trajectories of model (\ref{eq_eigant}), with $\alpha=\beta=0.1$, and $1 = d_1 = d_2 = d_3 > d_4 > \cdots > d_{10}$, illustrating Theorem \ref{thm_globstab}, equation (\ref{eq12}), namely, that the sum of states $x_1, x_2, x_3$ tends to the constant $(\beta d_1)/\alpha = 1$.}
\label{fig:illus_thm4_eq35}
\end{center}
\end{figure}

Finally, in order to illustrate (\ref{eq12}) of Theorem \ref{thm_globstab}, path lengths are chosen as $(1,1,1,2,3,4,\ldots,8)\in \R^{10}$, which means that the sum of the components $\sum_{k=1}^3 x_k (t)$, as $t$ tends to infinity, should tend to the fixed value $\beta d_1/\alpha$, which is $1$, in this case, since $d_1=1, \alpha=\beta$. This is confirmed by the simulation shown in figure \ref{fig:illus_thm4_eq35}.

\begin{figure}[h]
\begin{center}
\includegraphics[width=8.4cm]{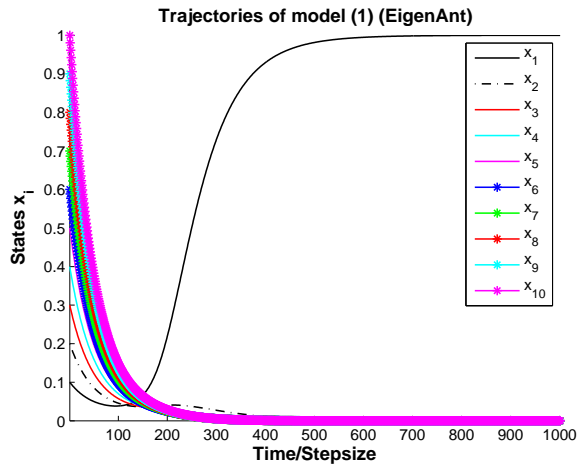}    % The printed column width is 8.4 cm.
\caption{Convergence to equilibrium corresponding to the shortest path, using EigenAnt dynamics (\ref{eq_eigant}), for a ten-path two-node shortest path problem described at the beginning of Section \ref{sec_simul}.}
\label{fig:eigenant_sync}
\end{center}
\end{figure}

%\begin{figure}[h]
%\begin{center}
%\includegraphics[width=8.4cm]{maxant_sync.eps}    % The printed column width is 8.4 cm.
%\caption{Convergence to equilibrium corresponding to the shortest path, using MaxAnt dynamics (\ref{eq_maxant}), for a ten-path two-node shortest path problem described at the beginning of Section \ref{sec_simul}.}
%\label{fig:maxant_sync}
%\end{center}
%\end{figure}

%\begin{figure}
%\begin{center}
%\includegraphics[width=8.4cm]{sigbant_sum_sync.eps}    % The printed column width is 8.4 cm.
%\caption{Convergence to equilibrium corresponding to the shortest path, using dynamics (\ref{}), for a ten-path two-node shortest path problem described at the beginning of Section \ref{sec_simul}.}
%\label{fig:sigbant_sum_sync}
%\end{center}
%\end{figure}

%Note that models (\ref{eq_phiant}), (\ref{eq_tanh_before}) have virtually indistinguishable behavior for the same choices of $\alpha,\beta,\phi$ and $\gamma (=10)$, whereas trajectories of model (\ref{eq_sign_before}) converge slightly faster than their counterparts (with the same choice of $\phi$).

\section{Concluding remarks}
This paper provided the first rigorous proof of global stability and asymptotic behavior of the continuous-time version of the discrete-time EigenAnt  dynamics. This is important because the discrete-time partially asynchronous version of the EigenAnt dynamics has been explored by \cite{JSBKC13} and shown to have essentially similar behavior (see \cite{JSBKC13} for details on the discrete-time partially asynchronous implementation). The implication is that the discrete-time analogs of the continuous-time variants proposed and studied in this paper, which can converge faster than EigenAnt, should also be useful for ACO applications, which are typically discrete-time and partially asynchronous. The property of robustness of stability with regard to parameter choices ($\alpha, \beta$) observed in \cite{JSBKC13} has been given a firm theoretical basis in the continuous-time case studied in the present paper. The EigenAnt algorithm was aptly referred to as a bare bones algorithm in \cite{EAW14}, which successfully incorporated it into the larger setting of ACO metaheuristics for solving multiple node shortest path problems such as the sequential ordering problem. An important reason for the bare bones terminology is that EigenAnt has only two parameters and, crucially, because of the global convergence property for all parameter choices, which was shown in the present paper, its performance is not critically dependent on these parameters. Thus, this paper has introduced a new class of bare bones algorithms that generalize EigenAnt and should therefore be of interest in a larger class of applications than the simple (paradigmatic) one that was subjected to a complete theoretical analysis herein.

\section*{Acknowledgment}
This research was supported by Inria (France), CNPq (Brazil), and the Dept.~of Science and Technology (India).

\end{document}